\documentclass[english]{article}
\usepackage[T1]{fontenc}
\usepackage[latin9]{inputenc}
\usepackage{esint}
\usepackage{babel}

\usepackage{float}
\usepackage{graphicx}
\usepackage{amssymb}

\begin{document}


\title{The Average Projected Area Theorem -- Generalization to Higher Dimensions}

\author{Zachary Slepian }

\maketitle

\nonumber \section{Abstract}
In 3-d the average projected area of a convex solid is $1/4$ the surface area, as Cauchy showed in the 19th century. In general, the ratio in $n$ dimensions may be obtained from Cauchy's surface area formula, which is in turn a special case of Kubota's theorem.  However, while these latter results are well-known to those working in integral geometry or the theory of convex bodies, the results are largely unknown to the physics community---so much so that even the 3-$d$ result is sometimes said to have first been proven by an astronomer in the early 20th century!   This is likely because the standard proofs in the mathematical literature are, by and large, couched in terms of concepts that are may not be familiar to many physicists.  Therefore, in this work, we  present a simple geometrical method of calculating the ratio of average projected area to surface area for convex bodies in arbitrary dimensions.  We focus on a pedagogical, physically intuitive treatment that it is hoped will be useful to those in the physics community.  We do discuss the mathematical background of the theorem as well, pointing those who may be interested to sources that offer the proofs that are standard in the fields of integral geometry and the theory of convex bodies.  We also provide discussion of the applications of the theorem, especially noting that higher-dimensional ratios may be of use for constructing observational tests of string theory.  Finally, we examine the limiting behavior of the ratio with the goal of offering intuition on its behavior by pointing out a suggestive connection with a well-known fact in statistics.

\section{Introduction}

	It is well known that for an arbitrary convex solid in three dimensions, the average projected
area is one-fourth the surface area.  Perhaps less well-known is the story of how this theorem came to be known in the physics community.  Karl Schwarzschild, a German astronomer primarily famous for his black hole solution to Einstein's general relativity, proved the theorem near the turn of the twentieth century.  Perhaps due to his untimely death at forty-two in 1916, however, the theorem was not well-known.  So it remained for his son, astronomer Martin Schwarzschild, to prove the theorem after the second World War, unaware of his father's work  (Knapp, personal communication).\footnote{The younger Schwarzschild is known in astronomy for his contributions to understanding stellar evolution, especially red giant stars.}  One can only imagine his surprise to learn from a colleague that his father had already proved it decades earlier!  

Meanwhile, the result in 2 and 3 dimensions had been well-known in the mathematics community since Cauchy's proof in works of 1841 and 1850.  The result in arbitrary dimension is known as ``Cauchy's surface area formula'' for this reason, although he did not prove it for higher dimensions.  But this was done by many others in the years following: Minkowski, Kubota, and Bonnesen all offer proofs (Bonnesen \& Fenchel 1987; see also Webster 1994). Essentially, the result may be seen as a special case of Kubota's theorem (see Klain \& Rota 1997, pp. 125-128; also Schneider 1993, pp. 295).  

Klain \& Rota present a so-called "mean projection formula" separately from Cauchy's surface area formula, and interpret the former as a form of Hadwiger's formula for compact convex sets (cf. pp. 94, 91, 125).  The proofs used involve objects such as Grassmannian measures and Haar probability measures, which may be unfamiliar to many physicists. The proofs in Bonnesen \& Fenchel (1987) and Eggleston (1958) are very similar to each other, and take a perhaps more elementary approach, but neither offers much intuition on why the formula holds.  Schneider (1993) simply proves a general mean value formula for projections, points out that what he terms Kubota's integral recursion (Kubota's theorem in Klain \& Rota) is a special case, and that a further specialization of the latter offers Cauchy's surface area formula.  Schneider's approach involves integrating over the group of rotations, as well as so-called mixed volumes and quermassintegrals, and we expect these may also be unfamiliar to many physicists. Santalo's (1976) discussion is perhaps the most geometric and physically intuitive of the standard references, and that closest to our approach here.  Again, though, certain elements of the proof that may not be obvious are asserted without further detail, and not much physical picture is provided to offer intuition as to why the formula holds.

Therefore, it is hoped that the present work, while not an original result, will be of use especially to those in the physics community who may have little background in integral geometry or the theory of convex bodies.  With this in mind, we spend some words developing a physical picture corresponding to the theorem in 3-$d$, and then generalize this picture to provide an intuitive proof in higher dimensions.  We also discuss applications of the theorem, both in 3-$d$ and in higher dimensions.  Finally, we provide a simple formula for calculating the ratio of average projected area to surface area in $n$ dimensions, and discuss the limiting behavior of this formula, both items unmentioned in the standard references discussed above.

\section{Definition of problem}

There are two equivalent, convenient ways to imagine average projected area. Consider an arbitrary convex solid in $n$ dimensions.  Now, hold the solid's position fixed, and shine a light beam
behind it toward an observer in front of it. Vary the light beam's
direction (and consequently the observer's location, since two points
define a line and so the observer must move to see the solid's shadow) through
all possibilities, measure the area of the shadow cast at each, and
average (Figure 1). Alternatively, one may hold the light beam (and observer's) position
fixed and rotate the solid through all possible orientations, keeping
its spatial location fixed. Measure the shadows and average (Figure 2). Either
method yields the same result -- as one would expect, since they
correspond to the same process but observed either in the rest frame
of the solid or of the observer.

\begin{figure}       
\includegraphics[scale=.4]{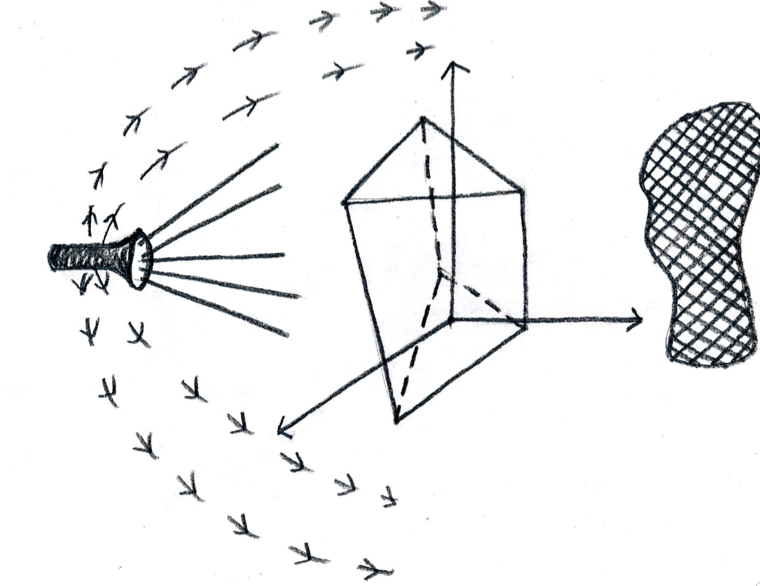} 
\centering{}\caption{One may imagine the average projected area as the average shadow cast by the solid when the observer's position is varied through all possible orientations with respect to the solid. The arrowed, dashed lines indicate the flashlight/observer's changing its position; the hatched shape at right is one example of the shadow that would be cast.}
 \vspace{0.5in}
\label{f:movingobserver}
\end{figure}

\begin{figure}       
\includegraphics[scale=.4]{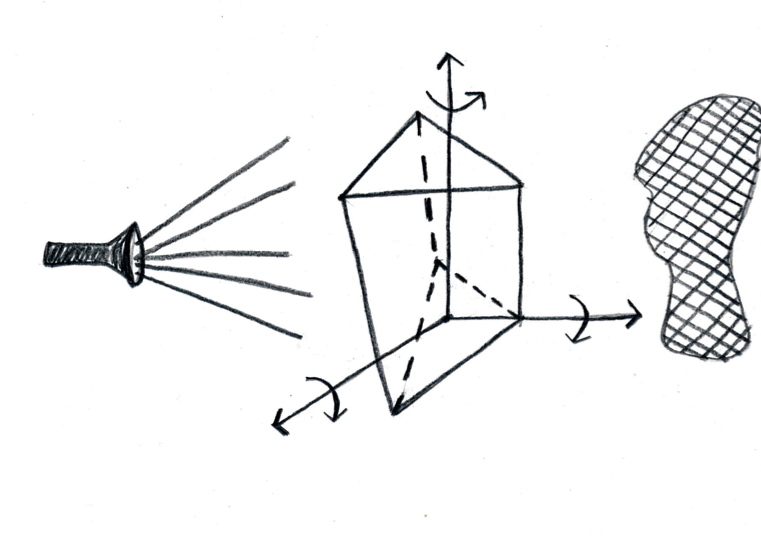} 
\centering{}\caption{One may also imagine the average projected area as the average shadow cast by the solid when it is rotated through all possible orientations with respect to a fixed observer. The arrowed semi-circles on the coordinate axes indicate the solid's rotation; the hatched shape at right is one example of the shadow that would be cast.}
 \vspace{0.5in}
\label{f:movingsolid}
\end{figure}

\section{Applications of the projected area theorem}
Applications of the theorem
are many. In astrophysics, dust is along the path light takes from
distant stars or galaxies to us, and so extinguishes some of the light.
The dust grains are irregularly shaped and oriented, so it is necessary
to correct for their average extinction of the incoming light. This
theorem provides the most convenient way to do so. 
The next step is complex numerical simulations; see e.g. Draine 2011.

The theorem also readily allows one to calculate the temperature of dust grains (or asteroids) illumined by a star (assuming the grains or asteroids are convex).  The energy received by a dust grain is just proportional to its area as projected onto the light source's: clearly, the grain cannot receive energy from any part of its surface not ``visible'' to the illumining star.  However, the particle can radiate energy away proportionally to its surface area.  Setting energy absorbed per unit time equal to energy radiated per unit time (i.e. the luminosities $L_{\rm abs}$ and $L_{\rm rad}$) determines the equilibrium temperature:
\begin{equation}
L_{\rm rad}=A_S \sigma_{\rm SB} T_{\rm grain}^4=L_{\rm abs}=\left<A_{\rm proj}\right>\sigma_{\rm SB} T_*^4\frac{4\pi R_*^2}{4\pi d^2}(1-a).
\end{equation}
$A_S$ is the surface area of the grain, $\sigma_{\rm SB}$ the Stefan-Boltzmann constant, $T_{\rm grain}$ the grain temperature, $\left<A_{\rm proj}\right>$ the average projected area of the grain, $T_*$ the star's temperature, $R_*$ the star's radius, $d$ the distance between the grain and the star, and $a$ the albedo, a measure of how reflective of light the grain is (the higher $a$, the less energy the grain absorbs).  Clearly, the ratio $\left<A_{\rm proj}\right>/A_S$, is required to calculate the temperature of the grain.

The average projected area theorem also has more terrestrial applications.  It is a starting point
for determining the 3-d shape of micro-particles using a 2-d map of
their projected areas (see Vickers \& Brown 2001, Brown et al. 2005).
The theorem has military applications as well: Saucier (2000) uses it
in discussion of simulating debris fragments behind armor and characterizing
the penetration potential of projectiles.  Further, it is a useful approximation
when considering particle transport in 3-d; see Glassner 1995. The
theorem also appears as a step in estimating the sphericity of fruits
and of polymer beads in gels; see e.g. Houston 1957 and Nussinovitch
2010. Finally, the theorem has applications in ray tracing for computer
graphics (Hanrahan 2000).

More speculatively, it is possible that measurements of the ratio of average projected area to surface area might offer insight on the number of dimensions in which we live.  String theory suggests there may be 9 (supersymmetric theories), 10 (M-theory), or 25 spatial dimensions (bosonic string theory); for reviews, see Kiritsis 1998, Polchinski 1998, or Schwarz 2000.  However, in string theory these extra dimensions
are often constrained to be on very small scales -- for instance, in Kaluza-Klein theory, the extra dimension must be smaller than $10^{-16}$ cm (Yagi et al. 2011), and in Randall-Sundrum models (4 spatial dimensions, see Randall \& Sundrum 1999 I and II), constraints are on the order of micro ($10^{-6}$) meters (Kapner et al. 2007). Constraining extra spatial dimensions using this theorem might also be challenging because in string theory, electromagnetic waves propagate through the usual 3 dimensions; it is only gravitational effects that are hypothesized to propagate in extra dimensions, as this would resolve the so-called hierarchy problem by explaining why gravity is very much weaker than electromagnetism.  Thus, it would be necessary to find applications in which the gravitational wave emission or absorption of an object depended on its average projected area.

\section{Plan of what follows}
The main purpose of this paper is to give an intuitive, geometric proof of the generalization of this theorem to higher-dimensional
convex solids and obtain an expression for $<A_{proj}>/A_S$ in arbitrary dimensions. A second goal will be to consider and offer intuition for the limiting behavior of this ratio as a function of dimension.

The plan of what follows is thus. I begin
by proving the theorem in 3-d for completeness, though this result may be found be found elsewhere, e.g. in Hildebrand 1942. We then present a geometric method for computing the ratio of average projected area to surface area as
a function of the dimension of the space (\S7), provide a table of values, and examine the limit of the ratio as the dimension becomes infinite (\S8). We conclude by presenting a simple recursion relation for the ratio between two consecutive dimensions and supplying an explicit formula for the ratio in terms of elementary operations.

\section{Proof in 3-d}

Consider on an arbitrary convex solid an infinitesimal unit
of surface area $dA$. Define $\vec{dA}$ as the unit vector normal to
the surface. Now consider a sphere of unit radius centered on the solid.  Consider an observer at a point on this sphere's surface, and let the observer's location define a hemisphere with its apex at their location. 
If the vector $\vec{dA}$ points such that it intersects this hemisphere, then the observer will see it. Since this hemisphere's
position is uncorrelated with the direction of $\vec{dA}$ for an
arbitrary unit of surface area $dA$, on average $\vec{dA}$ will
point into the observer's hemisphere half the time. 

Now restrict attention to the case where $\vec{dA}$ points
into the observer's hemisphere. Consider the projection of $\vec{dA}$ along the line
of sight from the surface to the observer. Denote the angle between $\vec{dA}$ and the line of sight
by $\theta,$ and note that because we have restricted ourselves
to a hemisphere, $\theta$ runs from zero to $\pi/2$.

We wish to integrate $\cos\theta$ to find its average value, as $\cos\theta$
gives the projection of $\vec{dA}$ onto the line of sight:

\begin{equation}
\left<\vec{dA}\cdot\hat{z}\right>=\int_{0}^{\pi/2}\cos\theta\sin\theta d\theta=\frac{1}{2}.
\label{e:3dfirstint}
\end{equation}

Notice that this integral is simply the differential projection of the solid's face onto the hemisphere at whose apex the observer sits, normalized by the hemisphere's surface area.

Here, the intuitive picture is that the orientation of the solid's face is changing
while the observer's position is held fixed. In this picture, we are in the observer's rest frame. As noted \S5,
it is equivalent to consider the solid's face fixed and vary the observer's
position with respect to it---i.e. to work in the solid's rest frame. The only quantity of physical significance is
the relative angle between the line of sight and the oriented unit
area vector $\vec{dA}$. 

Recalling  the additional factor of $1/2$ because $\vec{dA}$
only points into the observer's hemisphere half the time, and integrating over the solid's surface, we find

\begin{equation}
\left<A_{\rm proj,\:3-d}\right>=\int_{S}\left<\vec{dA}\cdot\hat{z}\right>dS=\frac{1}{4}A_{S},
\label{e:3dsecondeq}
\end{equation}

with $S$ denoting an integral over the solid's surface and $A_{S}$ the solid's surface area.

Note that, consonant with our comments above, though we derived equation (\ref{e:3dfirstint}) in
the rest frame of the observer, working in the rest frame of the solid
will prove more fruitful for generalizing the theorem to higher dimensions.
For notice that, in this picture, as the observer's position varies
(at fixed distance from the origin), it traces out a hemisphere. Hence we simply require the normalized projection of the solid's face onto
this hemisphere. But since the projection operator is commutative,
this is just the projection of the hemisphere onto the solid's
face.  This latter is simply a disc bounded by the circle produced where the hemisphere intersects the solids face.

\section{Proof of method for higher dimensions}

As noted in \S6, considering
the rest frame of the solid's face and projecting onto the hemisphere
defined by the different possible orientations of an observer at fixed
distance from the origin is a fruitful way to calculate the average projected area.
As also noted, this projection is easily evaluated by
finding the normalized projection of the hemisphere onto the solid's face. Both
points also hold in arbitrary higher dimensions.

In higher dimensions, one simply considers the hyperhemisphere generated
by varying the observer's orientation (at fixed distance from the origin) with respect
to the face of the solid under consideration, and calculates the normalized projection of the solid's hypersurface
onto it. This is equal to the normalized projection of the hyperhemisphere
onto the hypersurface bounding it. Now, the projection of a hyperhemisphere of dimension
$d$ onto the bounding hypersurface in $d$ dimensions is simply the hypersphere
of dimension $d-1$. Hence

\begin{equation}
\left<dA_{\rm proj,\: d}\right>=\frac{V_{d-1}}{2S_{{\rm H},\;d}},
\label{e:secproof}
\end{equation}

with $V_{d-1}$ the volume of the unit hypersphere of dimension $d-1$ and $S_{{\rm H},\;d}$ its surface area.  The factor of $1/2$ comes from the fact that $\vec{dA}$ only points into the observer's hemisphere half the time.

\section{Calculation in higher dimensions}

Here, we use the observation of the previous section to derive the
analog of equation (\ref{e:3dsecondeq}) for dimensions other than 3. As equation (\ref{e:secproof})
shows, the differential average projected area of a $d$-dimensional convex solid
is just the volume of the unit $d-1$-sphere normalized by double the surface area of the unit $d$-hemisphere. Hence we seek the factor
$k(d)$ relating $\left<A_{\rm proj, \;d}\right>$ to $A_S$, the surface area of the convex solid.  As we will see below, it is just $V_{d-1}/S_{d},$ where $V_{d-1}$ is the volume
of the unit $d-1$-sphere and $S_{d}$ is the surface area of the unit $d$-sphere.
We have

\begin{equation}
\left<A_{\rm proj, \;d}\right>=\int_{S}\left<dA_{\rm proj,\: d}\right>dS=\frac{V_{d-1}}{S_d}A_S=k(d)A_S,\end{equation}

where the second equality is from using equation (\ref{e:secproof}) for $\left<dA_{\rm proj,\: d}\right>$ and replacing $2S_{{\rm H},\;d}$ with $S_d$.
The result is easily verified for 3-d: with $d=3,$ $V_{d-1}=V_{2}=\pi$ and $S_{d}=4\pi$.

Now, for a sphere, 

\begin{equation}
V_{d-1}=\frac{S_{d-1}}{d-1},
\end{equation}

and

\begin{equation}
S_{d-1}=\frac{2\pi^{d/2-1/2}}{\Gamma\left(\frac{1}{2}d-\frac{1}{2}\right)}=\frac{1}{M_{d}\sqrt{\pi}}S_{d}\end{equation}

(Hypersphere, Wolfram), where we have defined $ $

\begin{equation}
M_{d}=\frac{\Gamma\left(\frac{1}{2}\left(d-1\right)\right)}{\Gamma\left(\frac{1}{2}d\right)}.\end{equation}

We can thus see that 

\begin{equation}
k(d)=\frac{1}{\sqrt{\pi}(d-1)M_{d}}.
\label{e:kd}
\end{equation}



Equation (\ref{e:kd}) is easily evaluated; we provide a table of the first thirty-two values, as well as a plot illustrating the monotonic decrease in $k(d)$.  For large $d$, $k\propto 1/\sqrt{d}$, as we show by a series expansion (equation (\ref{e:kdseries})).  This is a suggestive result, as it recalls the theorem in statistics that the standard deviation $\sigma$ of a set of measurements falls as $1/\sqrt{N}$ for $N$ measurements.  One may make this connection explicit by considering $N$ identical free particles moving in 3 spatial dimensions and examining the point set defined by the tip of each particle's velocity vector.  Connecting this point set results in a 3 dimensional convex solid.  One may now consider the area of a face of this solid.  From dimensional analysis, it must be proportional to the square of the separation between the vectors corresponding to the particles' whose coordinates define the face.  This separation is simply $\sigma$, with $\sigma$ the standard deviation of the velocity distribution from which the particles were drawn.  Thus the area of a face is proportional to $\sigma^2$. Now, there will be of order $N$ faces of the solid, so by equipartition arguments, on average the area of one face will be a fraction $1/N$ of the surface area.  Therefore we see that $\sigma \propto 1/\sqrt{N}$. This gives some intuition for the connection between the fact in statistics that $\sigma\propto 1/\sqrt{N}$ and the fact that, in many dimensions $d$, which one may interpret as describing a many ($N$) particle system with $d=N$, the ratio of average projected area to surface area also falls with $1/\sqrt{d}$.

\begin{table}[!ht] 
 \vspace{0.2in}
\caption{Ratio of average projected area to surface area ($k(d)$) as a function of dimension $d$}
 \vspace{0.04in}
\centering
    \begin{tabular}{ |p{.4 cm}|p{1.2 cm} ||p{.4 cm} |p{1.2 cm}|}
 
    \hline
    $d$ & $k(d)$ & $d$ & $k(d)$ \\ \hline
     2 &.318 &18 & .095  \\ \hline
     3 & .250 &  19 & .093  \\ \hline
     4 & .212  &20& .090  \\ \hline
     5 & .1875  &21& .088  \\ \hline
     6 & .170 &22& .086  \\ \hline
     7 & .156 & 23& .084  \\ \hline
     8 & .146  &24& .082  \\ \hline
     9 & .137  & 25& .081  \\ \hline
     10 & .129 &26& .079  \\ \hline
     11 & .123  &27& .077  \\ \hline
     12& .118&28& .076  \\ \hline
     13& .113 &29& .075  \\ \hline
     14 & .109 &   30& .073  \\ \hline
     15 & .105  &31& .072  \\ \hline
     16 & .101  &32& .071  \\ \hline
     17 & .098  &   33& .070  \\ \hline

      \end{tabular}
     \vspace{0.2in}
     \label{t:kd}
\end{table}

\begin{figure}
\includegraphics[scale=1.2]{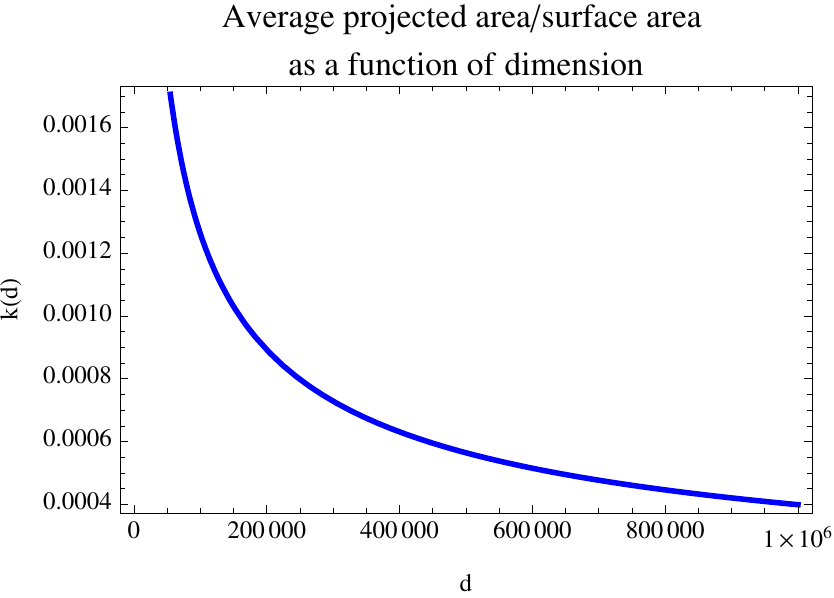}
\centering{}
\caption{ This plot shows the monotonic decrease in  $k(d)$ with increasing $d$; to leading order $k(d)\propto1/\sqrt{d}$, as equation (\ref{e:kdseries}) indicates.}
 \vspace{0.5in}

\label{f:kdbig}
\end{figure}

Finally, it is interesting to consider the limiting behavior of $k(d)$ as $d$ tends to infinity.  As might be intuited from the plot, $k(d)\rightarrow 0$ as $d\rightarrow \infty$.  This follows immediately from the series for $k(d)$ about $d=\infty$, with $x\equiv 1/d$:
\begin{equation}
k(d)=\frac{1}{\sqrt{2\pi}}\left\{x^{1/2}+\frac{x^{3/2}}{4}+\frac{x^{5/2}}{32}-\frac{5x^{7/2}}{128}-\frac{21x^{9/2}}{2048}+\mathcal{O}\left(x^{11/2}\right)\right\}.
\label{e:kdseries}
\end{equation}
This series is accurate to one part in ten thousand for $d$ as low as five.

\section{Recursion relation and explicit formula for $k(d)$}

There exists a simple recursion relation for $k(d+1)$ in terms of
$k(d)$. Observe that

\begin{equation}
M_{d}M_{d+1}=\frac{\Gamma\left(y\right)}{\Gamma\left(y+1\right)},\end{equation}

where $y=\frac{1}{2}(d-1)$ and $M_d$ is defined by equation (8). Using the identity that $\Gamma(1+y)=y\Gamma(y)$,
we have

\begin{equation}
M_{d}M_{d+1}=\frac{2}{d-1}.\end{equation}

Using equation (\ref{e:kd}) in the above to relate $k(d)$ and $k(d+1)$
to, respectively, $M(d)$ and $M(d+1),$ we find

\begin{equation}
k(d+1)=\frac{1}{2\pi dk(d)}.
\label{e:recursion}
\end{equation}

It is easily verified that, beginning with $k(3)=1/4$, this formula reproduces the results of equation (\ref{e:kd}) as given in Table 1. 

Using the recursion (equation (\ref{e:recursion})) and anchoring it at $k(2)=1/\pi$
allows derivation of an explicit formula for $k(d),$ $d>2$: 
\begin{equation}
k(d)=\frac{1}{2}\prod_{n=0}^{(d-3)/2}\frac{2n+1}{2n+2},\ \ d\ odd\end{equation}

and

\begin{equation}
k(d)=\frac{1}{\pi}\prod_{n=0}^{(d-4)/2}\frac{2n+2}{2n+3},\ \ d\ even.\end{equation}

\section{Acknowledgements}

I thank G. Knapp, J.R. Gott, and J. Pardon for useful discussion during the course of this work, and D. Yang and M. Reitzner for useful correspondence.

\section{References}

$\;\;\;\;\;$

Brown D.J., Vickers G.T., Collier A.P., \& Reynolds G.K., 2005, Chemical
Engineering Science, vol. 60, issue 1, p. 289.

Cauchy, A., 1841, Note sur divers theoremes relatifs a la rectification des courbes et a la quadrature des surfaces, C.R. Acad. Sci., Paris, vol. 13, pp. 1060-1065.

Cauchy, A., 1850, Memoire sur la rectification des courbes et la quadrature des surfaces courbes. Mem. Acad. Sci. Paris vol. 22, 1850, pp 3ff.

Draine B., 2011, Physics of the Interstellar and Intergalactic Medium.
Princeton University Press, Princeton. 

Eggleston H., 1958, Convexity. Cambridge University Press, New York.

Gamma Function, 2011, Wolfram MathWorld, 

http://mathworld.wolfram.com/GammaFunction.html

Gardner R.J., 2006, Geometric Tomography (2nd Ed.). Cambridge University Press, New York.

Glassner A.S., 1995, Principles of digital image synthesis, vol. 2,
p. 601.

Gruber P., 2007, Convex and Discrete Geometry. Springer Verlag, New York.

Hanrahan P., 2000, CS348B Lecture 8, http://graphics.stanford.edu/courses/cs348b-00/lectures/lecture03/raytracingII.pdf.gz.

Hildebrand R.H., 1983, Royal Astron. Soc. Quart. Jrn. V.24, No. 3,
p. 267, http://articles.adsabs.harvard.edu//full/1983QJRAS..24..267H/0000282.000.html

Houston R.K., 1957, Agric. Eng. 39: 856-858.

Hypersphere, 2011, Wolfram MathWorld, http://mathworld.wolfram.com/Hypersphere.html

Kapner D. J., Cook T. S., Adelberger E. G., Gundlach J. H.,
Heckel B. R., Hoyle C. D., \& Swanson H. E., 2007, Phys. Rev.
Lett. 98, 021101

Kiritsis E., 1998, arXiv: hep-th/9709062v2

Klain D. \& Rota G.-C., 1997, Introduction to Geometric Probability, Cambridge University Press, New York.

Knapp, G.R., 2011, personal correspondence

Nussinovitch A., 2010, Polymer Macro- and Micro-Gel Beads: Fundamentals
and Applications. Springer, New York, p. 15.

Polchinski J.,  1998, ``String theory. Vol. 1 and 2,'' Cambridge, UK: Univ. Pr.

Projection, 2011, Wolfram MathWorld, http://mathworld.wolfram.com/Projection.html

Randall L. \& Sundrum R., 1999, Phys. Rev. Lett. 83, 3370. (I)

Randall L. \& Sundrum R., 1999, Phys. Rev. Lett. 83, 4690. (II)

Santalo L., 1976, Integral Geometry and Geometric Probability. Addison-Wesley, Reading, MA.

Saucier R., 2000, Shape Factor of a Randomly Oriented Cylinder, ARL-TR-2269.

Schneider R., 1993, Convex Bodies: The Brunn-Minkowski Theory, Cambridge University Press, New York.

Schwarz J.H., 2000, arXiv: hep-ex/0008017

Vickers G.T. \& Brown D..J., 2001, Proc. R. Soc. Lond. A, 457, 283-306.

Webster R., 1994, Convexity. Oxford University Press, New York.

Yagi K., Tanahashi N., \& Tanaka T., 2011, Phys.Rev.D83: 084036.

\end{document}